\documentclass{article}
\usepackage[utf8]{inputenc}
\usepackage{graphicx}
\usepackage[english,russian]{babel}

\newcommand{\trn}{^{\scriptscriptstyle {\rm T}}}

\title{УПРАВЛЕНИЕ ВЫРАВНИВАНИЕМ ИНВАРИАНТОВ
\author{А.Л. Фрадков\\
Институт проблем машиноведения РАН, \\
61, Bolshoy pr .V.O., Saint Petersburg 199178 Russia \\
St Petersburg University \\ 28 Universitetskii prospect, Peterhof,
St Petersburg 195904, Russia\\
E-mail: fradkov@mail.ru
}
}
\date{Nov 24, 2023}

\begin{document}
\maketitle



Обсуждается новый класс задач управления - управление гомеостазисом.  Задачи управления гомеостазисом можно рассматривать как задачи управления с заданным целевым множеством, в частности, как задачу стабилизации значений некоторой целевой функции, являющейся инвариантом неуправляемой системы. В докладе рассматривается более общий класс: выравнивание значений двух или более целевых функций, каждая из которых является инвариантом соответствующей подсистемы сложной системы. Предложен
подход к синтезу алгоритмов управления на основе метода скоростного градиента
и установлены условия применимости подхода.



\section{Введение}

Гомеостазис (гомеостаз) – удивительное свойство живых существ и сообществ, которое остается не до конца понятым и по сей день.
Впервые термин {\it гомеостазис} ввёл в 1929 году американский физиолог Уолтер Кэннон (от греч. «хомойос»- тот же, подобный, и «стазис»- состояние). Мысль о том, что постоянство внутренней среды является необходимым условием для нормального существования любого живого организма высказал ещё в 1857 году французский физиолог Клод Бернар. Понятие гомеостазиса играет важную роль в биологии, медицине и других науках и означает поддержание
состояния внутренней среды системы в области нормального функционирования при изменении внешних условий. Если гомеостазис в системе не может быть обеспечен внутренними ресурсами, возникает вопрос о обеспечении его за счет внешних сил - управлении гомеостазисом.

Существующие в организмах принципы регуляции биологических процессов имеют много общего с принципами управления в неживых, технических системах. И в том, и в другом случае стабильность системы достигается с помощью определенной формы управления. Поэтому задача об управлении гомеостазисом относится к области кибернетики – науки об управлении и связи в живом организме и машине, как ее определил Н.Винер в 1948 г \cite{Wiener48}.

Понятия гомеостазиса и гомеостата приобрели популярность в кибернетике в 1950-х гг, когда английский психиатр и кибернетик У.Р.Эшби построил машину, демонстрировавшую устойчивое функционирование при изменении внешней среды и связей между элементами  и назвал ее гомеостатом \cite{Ashby49,Ashby60}. Как отмечал Г.Паск \cite{Pask61}, гомеостат Эшби, в отличие от многих имитаций, внешне ведущих себя подобно мозгу, был создан, чтобы понять внутреннее устройство мозга. Поэтому  гомеостат Эшби до сих пор привлекает внимание специалистов  и считается одним из прототипов адаптивного регулирования в природе и технике \cite{Herrmann04}. Собственно, о связи гомеостата с адаптивным поведением писал еще сам Эшби \cite{Ashby60}. Однако развитие адаптивного управления в XX веке позволило рассмотреть новые модели и методы, устанавливающие более многообразные связи управления и гомеостазиса.

 В \cite{Herrmann04} указано, что первые математические модели динамики гомеостата представляли собой линейные дифференциальные системы $dx/dt= A(k)x$ с коэффициентами $k$, переключающимися при достижении пороговых значений углов поворота магнитов и останавливающимися, когда система становится асимптотически устойчивой $Re \lambda_i(A(k))<0.$ В качестве альтернативы в \cite{Herrmann04} предлагалось использовать непрерывные алгоритмы адаптивного управления. Интересно, что алгоритмы \cite{Herrmann04} хорошо известны в теории адаптивного управления и их свойства могут быть установлены на основе метода функций Ляпунова и пассификации \cite{F76,ASF88}.

Во второй половине XX века процессы гомеостазиса интенсивно изучались в медицине, при этом активно использовались математические модели. Построение математических моделей позволяет связать представления биологов, математиков и инженеров и приблизиться к практическому использованию теоретических исследований в медицине. Было разработано множество математических моделей гомеостатического поведения конкретных подсистем организмов. Среди стабилизируемых показателей чаще всего фигурируют уровни глюкозы и инсулина в крови \cite{Lombarte13, Gaohua09}, уровни кальция \cite{Raposo02}, межклеточного железа \cite{Chifman12} и ряд других показателей.
Модели гомеостазиса конкретных биологических процессов обычно достаточно сложны, поскольку учитывают специфику нескольких взаимосвязанных процессов. Например, модель гомеостазиса кальция \cite{Raposo02} включает 11 нелинейных дифференциальных и ряд алгебраических уравнений. Математическая модель гомеостазиса межклеточного железа содержит 5 нелинейных дифференциальных упавнений и 15 неопределенных параметров \cite{Chifman12}.  Математическая модель дифференцировки $T$-клеток содержит 18 параметров \cite{Thomas08} и т.д.

Явление гомеостазиса может пониматься настолько широко и имеет настолько важные методологические аспекты, что в последние несколько десятилетий появилось множество работ, обобщающих представления о гомеостазисе на другие классы систем: технические, общественные и др.
В  целом ряде исследований  гомеостатические принципы из иммунологии применялись к проектированию роботов \cite{Vargas05,Kingson19}, сетей связи \cite{Oka14}  и других технических систем.
В работе \cite{Samigulina22} было предложено применять идеи гомеостазиса для разработки интеллектуальной технологии управления сложными технологическими объектами в нефтегазовой отрасли, основанной на единой искусственной иммунной системе (ЕАИС) и принципах иммунологического гомеостаза с учетом современной микропроцессорной техники, широко применяемой в современном индустриальном производстве.
Наиболее общие представлении о математических  моделях гомеостазиса были развиты недавно в прикладной математике. В работе \cite{Golubitsky17} дается следующее общее определение.

{\it Гомеостазис имеет место в биологической или химической системе, когда некоторая выходная переменная остается приблизительно постоянной при изменении входной переменной  на некотором интервале}.

Разумеется, это определение применимо не только к биологическим или химическим системам. В работе \cite{Golubitsky17} предлагается уточнение этого определения. Отмечается, что интервал изменения входных параметров (внешних условий), в котором некоторая выходная переменная $y$ остается постоянной или  приблизительно постоянной, может быть достаточно большим, т.е. $y$ может являться инвариантом или приближенным инвариантом системы относительно внешних воздействий. Перспективу еще большего расширения областей инвариантности при гомеостазисе предоставляет введение управляющих переменных и применение методов теории управления. Дальнейшее обобщение может быть направлено на поддержание при гомеостазисе не значения какого-то инварианта, а некоторого соотношения между несколькими инвариантами системы.

В настоящей работе рассматривается класс систем, обладающих инвариантами и предлагается общий способ управления, обеспечивающий сохранение соотношений между инвариантами системы  при изменении начальных условий и параметров внешней среды.


\section{Постановка задачи}

В общем случае под гомеостазисом будем понимать способность организмов и систем сохранять свое состояние в области устойчивого функционирования при изменении внешних условий. Уточняя сказанное выше, обозначим набор переменных, характеризующих состояние системы в момент времени $t$, через $x(t)\in R^n$, а область устойчивого функционирования обозначим через $X^*\subset R^n$. Тогда условие гомеостаза можно записать в виде
\begin{equation} \label{1}
x(t)\in X^*.
\end{equation}
Рассмотрим систему описываемую  дифференциальными уравнениями состояния
\begin{equation} \label{2}
dx/dt= F(x(t),v(t)), ~y(t)=h(x(t))
\end{equation}
где $v(t)\in D\subset R^d$ - вектор внешних воздействий (параметров  системы или внешней среды),
  $y(t)\in R^l$ вектор измеряемых выходов системы.

  Пусть для любого $v\in D $ существует вектор $\bar x(v)$, являющийся равновесием
   системы $F(\bar x(v),v)=0$. Будем говорить, что система обладает гомеостазом
(сильным гомеостазом) по $y$ в $D$,
 если $\nabla_v y(\bar x(v))=0$ для всех $v\in D $.
Иначе говоря, это означает, что функция $y(\bar x(v))$ явлется инвариантом (сохраняет свое значение)
 при изменении $v\in D $. Такое определение близко к формализации \cite{Golubitsky17}.

Теперь предположим, что на систему (\ref{2}) может воздействовать управление,
т.е. математическая модель системы имеет вид
\begin{equation} \label{3}
dx/dt= F(x(t),u(t),v(t)), ~y(t)=h(x(t)),
\end{equation}
где $u(t)\in R^m $ - вектор управляющих переменных. Пусть для любого $v\in D $
существуют векторы $\bar u(v), \bar x(v)$ такие, что
\begin{equation} \label{4}
 F(\bar x(v),\bar u(v),v))=0, ~\nabla_v y(\bar x(v))=0,
\end{equation}
при $v\in D $. Свойство (\ref{4}) означает, что $y(x(v))$ является инвариантом системы при воздействии управления.

Очевидно, проблема состоит в том, чтобы найти функцию управления $u^*=U^*(v)$, обеспечивающую
инвариантность выхода $y(x(v))$. Сложность проблемы в том, что значение $v$ обычно неизвестно
и найти правильное $u^*$ в каждый момент времени не представляется возможным. Поэтому цель управления заменяется на асимптотическую цель
\begin{equation} \label{5}
\lim_{t\to\infty}y(x(t))= y^*,
\end{equation}
для любых  $v\in D $, где $D$ - известное множество возможных значений $v$,
а управление $u(t)$ строится на основе измерения в каждый момент времени текущих значений
вектора состояния $x(t)$ или вектора выходов $y(t)$ (в более общей ситуации при выборе $u(t)$
могут использоваться также значения $x(s)$ или   $y(s)$ в прошлые моменты времени $s<t$).

Таким образом, мы пришли к стандартной задаче теории автоматического управления.
Ей посвящено огромное количество работ и она адекватна многим задачам достижения гомеостазиса. Однако общность задачи, с одной стороны, не дает возможности сформулировать ее конструктивное решение, а, с другой стороны, не учитывает специфику многих конкретных задач. Множество $X^*$ может вырождаться в точку, может задавать какие-то желаемые значения для переменных организма, определяющих состояние гомеостазиса и т.д. Нас интересует более общая постановка задачи, когда  гомеостазис рассматривается не на уровне отдельных переменных, а на уровне системы, и целью управления является обеспечение каких-то соотношений между переменными, гарантирующих работоспособность целой, возможно, сложной  системы.

Поэтому ниже рассматривается  модифицированная, более общая формулировка цели (\ref{5})  для класса моделей (\ref{3}). Она требует не стремления траектории системы в некоторую окрестность  какой-то желаемой точки, а выполнения определенного соотношения между переменными системы. Это соотношение можно задать в виде
\begin{equation} \label{6}
h(x(t))=y^*,
\end{equation}
где $h(x)$  - некоторая заданная функция состояния системы. Соотношение (\ref{3}) означает, что траектория системы принадлежит некоторой заданной поверхности в пространстве состояний системы. Такая задача также рассматривалась во многих работах, она соответствует стабилизации решений на заданной поверхности и может быть названа частичной стабилизацией или стабилизации по выходу. Ниже при некоторых предположениях предлагается  решение задачи, основанное на методе скоростного градиента \cite{AndrFr21}.

 Основное накладываемое предположение – существование у каждой подсистемы инварианта: функции $y_i=h_i(x)$, которая сохраняется, если внешние воздействия на подсистему не меняются. Существование инвариантов возможно не только у биологических систем, но и у многих физических, технических, социальных и др. систем. Типичный пример – механическая система, состоящая из нескольких консервативных подсистем: у каждой подсистемы имеется функция энергии, которая сохраняется, если внешние воздействия на подсистему не меняются.

  Пусть сложная система описывается совокупностью уравнений состояния
\begin{equation} \label{2a}
dx_i/dt= F_i(x_i(t),u_i(t),v_i(t)),
\end{equation}
где $x_i(t)\in R^{n_i}$ - векторы состояния $i$-й подсистемы, $u_i(t)\in R^{m_{i}}$ - векторы управляющих воздействий, $v_i(t)\in D_i\subset R^{l_i}$ - векторы возмущений (параметров внешней среды)$i=1,\ldots,N$.
Для таких систем естественно поставить в качестве цели управления согласованную работу всех подсистем сложной системы или выравнивание значений важнейших переменных, чтобы обеспечить сравнимые условия функционирования подсистем. Такая цель формулируется в виде цепочки равенств
\begin{equation} \label{3a}
h_1(x_1(t))=h_2(x_2(t))=\ldots=h_N(x_N(t)),
\end{equation}
где $h_i(x_i)$ -- функция выходов $i$-й подсистемы, $i=1,2,\ldots,N$,  $N$ -- число подсистем. Цель (\ref{3a}) можно обобщать еще дальше, задав ее как некоторое соотношение между выходными переменными, учитывающее, например, наличие  масштабных коэффициентов и различие размерностей переменных, но математической общности цель (\ref{2a}) не умаляет. Наконец, будем рассматривать асимптотическую цель
\begin{equation} \label{3b}
h_i(x_i(t))-h_{i+1}(x_{i+1}(t))\to 0, i=1=\ldots=N,
\end{equation}
где при $i=N+1$ считается, что  $i=1$.


\subsection{Управление инвариантами нелинейных аффинных систем}


Рассмотрим задачу управления, считая, что управляемая система аффинна по
входу (по управлению).
\begin{equation} \label{7}
\dot x=f(x)+g(x)u, ~y=h(x),
\end{equation}
где $x\in R^n, u\in R^m, y\in R^l$, а $f,g,h$ -- гладкие функции.
Введем обозначения $L_f h(x)=$, $Z(x)=\nabla h^{\trn}$
и предположим, что существует неотрицательная функция $V(x)\ge 0$, такая, что выполнены соотношения \cite{ShirFr00}.

(A1) Для любой точки $x$ фазового пространства системы $X\subset R^n$ выполняется неравенство $L_f(x)V(x)\leq 0$, означающее, что система устойчива по отношению к функции Ляпунова  $V(x)$.

(A2) Для любой точки $x$, принадлежащей множеству
$x\in X: L_g V(x)=0$, выполняется тождество $L_f h(x)\equiv 0$, 
означающее, что аффинная часть системы (\ref{7}) не допускает потери
на множесьве $\{x\in X^ L_g V(x)\equiv 0\}$.


(A3) $dim S(x)=l$ при $Z(x)=0, ~x\in \Omega_0$, где
$$
S(x)=span\{Z(x), L_f Z(x), L^2_f Z(x),\ldots,\}
$$

(A4) Существует $\varepsilon>0$, такое, что любая связная компонента множества $D_{\varepsilon}$ ограничена, где
$$
D_{\varepsilon}=\Omega_0\cap\{x:|det Z(x)Z(x)^{\trn}|<\varepsilon\}.
$$

Справедлив следующий результат.

{\bf Теорема 2.1 \cite{ShirFr00}.}
Пусть  выполнены условия (A1)-- (A4)
Then the goal $\lim_{t\to\infty} Q(x(t))=0$ is achieved for any $x(0)\in \Omega_0$.

{\bf Замечание 1}. Если $\Omega_0$ iограничено, то условие (A4) можно опустить.


\section{Управление соотношением инвариантов \\
 нелинейных систем}

Рассмотрим набор нелинейных, аффинных по управлению систем
 \begin{equation}
    \label{S1}
    \dot x_i=f_i(x_i)+g_i(x_i)u_i, y(i)=h_i(x_i), i=1,\ldots,N.
    \end{equation}

Предположим, что существует неотрицательная функция $V(x)\ge 0$, такая, что выполнены соотношения



(AN1) для любой точки  $x_i\in X\subset R^n$
выполняются неравенства
$L_{f_i(x_i)}V(x)\leq 0$, означающие устойчивость невозмущенных систем
по отношению к функции $V(x)$ (пассивность).

(AN2) для любой точки $x_i$, из множества
$x\in X: L_g V(x)=0$, такой, что тождества $L_g(x_i)V(x)=0$ выполняются при $L_{f_i(x_i)} h_i(x_i)\equiv 0$/

Условия (AN1), (AN2) являются обобщением условий (A1),(A2) на сложные системы.
Из них следует, что функция $h_i(x)$ является инвариантом $i$--й подсистемы.

Поставим задачу достижения при $t\to\infty$ заданного соотношения между выходами подсистем.
Не умаляя общности, ее можно сформулировать как обеспечение
    выполнения асимптотических соотношений
 \begin{equation}
    \label{S2a}
    y_i(t)-y_{i+1}(t)\to 0, ~~i=1,\ldots,N-1
    \end{equation}
при $t\to\infty$.  Эта задача является обобщением задачи управления инвариантами нелинейных систем, рассмотренной в предыдущем разделе.

Для решения задачи  применим метод скоростного градиента.
 Выберем целевую функцию в виде циклической суммы квадратов ошибок
 \begin{equation}
    \label{S3}
    Q(x)=(y_1-y_2)^2+(y_2-y_3)^2+\ldots (y_{n}-y_{n+1})^2,
    \end{equation}
    где считается, что $y_{n+1}=y_1$
    Очевидно, при $Q(x)=0$ все выходы подсистем совпадают и цель (\ref{S2a}) достигается

   Для синтеза управления применим модифицированный метод скоростного градиента,
   вычисляя градиент по $u$ не от скорости изменения целевой функции
   $\dot Q(x,u)$, а от ее верхней оценки:

  \begin{equation}
    \label{S8}
    u_i=-\gamma(2g_i\nabla h_i^{\trn})(2y_i-y_{i-1}-y_{i+1}), i=1,\ldots,N
      \end{equation}
Справедливы следующие условия достижения цели в синтезированной системе.

 {bf Теорема 1.} При выполнении условий (AN1), (AN2) для   системы  (\ref{S1}) с алгоритмом управления (\ref{S8}) выполняется неравенство $Q(x(t)) \leq Q(x(0))$, $\forall t \geq 0$ и обеспечивается соотношение $u(t) \to 0 $ при $t \to \infty$. Кроме того, обеспечивается альтернатива: на траектории $x(t)$ либо достигается цель управления (\ref{S2a}), либо 
  $g_i\nabla h_i^{\trn}\to 0$ при $t \to \infty$ для некоторого $i=1,\ldots,N$. Если, кроме того, выполнены условия (A3),(A4), то исходная цель
  (\ref{S2a})  достигается.

{\bf Доказательство}.
При использовании закона управления (\ref{S8}) и при выполнении условий (AN1), (AN2) выполнены соотношения
$$
\dot Q(x,u)=-\gamma\sum_{i=1}^N((2g_i\nabla h_i^{\trn})^2(2y_i-y_{i-1}-y_{i+1})^2\le 0,
$$
откуда функция $Q(x(t)$ ограничена: $Q(x(t)\le Q(x(0)$. Поскольку $Q(x(t))$
не возрастает, существует предел $\lim_{t\to\infty}Q(x(t))=Q_{\infty}$.
Если $Q_{\infty}=0$, то цель достигается и теорема доказана.
Если $Q_{\infty}> 0$, то, значит, хотя бы одно слагаемое вида
$(2y_i-y_{i-1}-y_{i+1})$ не стремится у нулю и, следовательно,
хотя бы одна функция $(2g_i\nabla h_i^{\trn})$ стремится к нулю.
При дополнительном выполнении условий (A3),(A4), аналогично \cite{ShirFr00}
устанавливается достижение исходной цели (\ref{S2a}).

\section{Заключение}
\hspace*{1.25cm}Исследована возможность применения метода скоростного градиента к задаче обеспечения заданного соотношения инвариантов в нелинейных системах.  Поставлена задача выравнивания уровней энергий в системе, получен теоретический результат о возможности достижения цели управления для произвольного числа подсистем. Полученные результаты могут иметь применения в медицине и в технике, поскольку  применимы как к живым организмам, так и к техническим системам.

В дальнейшем представляет  интерес исследование сетевой модели взаимодействия подсистем с наличием физических связей (пружин) между ними. Сложность в этом случае состоит в том, что каждая отдельная подсистема при нулевом управлении уже не является гамильтоновой системой, что приводит к необходимости оценки влияния взаимосвязей.

Работа поддержана грантом СПбГУ ID 94034465.


\end{document}